\documentclass[final,hidelinks,onefignum,onetabnum]{siamart251216}
\usepackage{lipsum}
\usepackage{amsfonts}
\usepackage{bookmark}
\usepackage{bm}
\usepackage{graphicx}
\usepackage{epstopdf}
\usepackage{hyperref}
\usepackage{enumitem}
\usepackage{algorithmic}
\usepackage{booktabs}
\usepackage{multirow}
\usepackage{cleveref}
\usepackage{aliascnt}
\ifpdf
  \DeclareGraphicsExtensions{.eps,.pdf,.png,.jpg}
\else
  \DeclareGraphicsExtensions{.eps}
\fi

\newsiamremark{remark}{Remark}
\newsiamremark{hypothesis}{Hypothesis}
\crefname{hypothesis}{Hypothesis}{Hypotheses}
\newsiamthm{claim}{Claim}
\newsiamremark{fact}{Fact}
\crefname{fact}{Fact}{Facts}

\newaliascnt{example}{theorem}
\newtheorem{example}[example]{\textit{Example}}
\aliascntresetthe{example}
\crefname{example}{Example}{Examples}
\Crefname{example}{Example}{Examples}

\crefname{enumi}{part}{parts}
\Crefname{enumi}{Part}{Parts}

\newcommand{\sd}{\mathord{\prime\mkern-2.5mu\reflectbox{$\scriptstyle\prime$}}}

\Crefname{ALC@unique}{Line}{Lines}

\headers{The specular gradient method}{K. Jung and J. Oh}

\title{Nonsmooth convex optimization using the specular gradient method with Root-linear convergence\thanks{Date: \today
\funding{This work was supported by the National Research Foundation of Korea(NRF) grant funded by the Korea government(MSIT) (Nos. NRF-2020R1C1C1A01014904, NRF-RS-2023-00217116).
This work is supported in part by Michigan State University and the National Science Foundation Research Traineeship Program (DGE-2152014) to Kiyuob Jung.}}
}

\author{
  Kiyuob Jung\thanks{Department of Mathematics, Michigan State University, East Lansing, MI 48824 USA (\email{kyjung@msu.edu}).}
  \and Jehan Oh\thanks{Department of Mathematics, Kyungpook National University, Daegu, 41566, Republic of Korea (\email{jehan.oh@knu.ac.kr}).}
}

\usepackage{amsopn}

\usepackage[normalem]{ulem}
\usepackage{placeins}

\definecolor{SPEGpurple}{RGB}{128,0,128}

\hypersetup{colorlinks=false, pdfborder={0 0 1}, linkbordercolor={1 0 0}, citebordercolor={0 1 0}}

\begin{document}
\maketitle

\begin{abstract}
  We propose the specular gradient method for one-dimensional convex optimization.
  Assuming that the minimum is attained and a suitable initial distance bound holds, we establish R-linear convergence using normalized steps of geometrically decreasing length.
  Neither strong convexity nor differentiability is required.
\end{abstract}



\section{Introduction}

Consider the one-dimensional convex optimization problem:
\begin{equation}
  \label{opt:1d}
  \min_{x \in \mathbb{R}} f(x),
\end{equation}
where $f : \mathbb{R} \to \mathbb{R}$ is a given convex objective function.
Since $f$ need not be differentiable, classical gradient methods may not apply to \eqref{opt:1d}.
For such nonsmooth problems, the subgradient method generates a sequence from a given initial point $x_0 \in \mathbb{R}$ by
\begin{displaymath}
  x_{k+1} = x_k - \gamma_k \, g_k, \qquad k = 0, 1, 2, \ldots,
\end{displaymath}
where $g_k \in \partial f(x_k)$ is a subgradient of $f$ at $x_k$ and $\gamma_k > 0$ is the step size.
Its convergence depends on the step sizes and the information available about the objective.
For background on subgradients and their applications to optimization, we refer to \cite{2014_Boyd,2022_Boyd,1977_Goffin,2018_Nesterov_BOOK,1985_Shor_BOOK}.

Related developments in nonsmooth optimization include optimality and duality results based on symmetric derivatives and generalized convexity \cite{1971_Minch}, and algorithms using $\varepsilon$-subgradients \cite{1973_Bertsekas}.
More general formulations use norms and duality pairings in place of the standard Euclidean geometry, with convex analysis extending to linear spaces \cite{2018_Nesterov_BOOK,1984_VanTiel_BOOK}.
For a broad treatment of convex optimization, see \cite{2004_Boyd_BOOK}.

Under suitable step-size and boundedness assumptions, standard subgradient methods have sublinear guarantees for the best objective value \cite{2014_Boyd,2018_Nesterov_BOOK}.
Faster rates are available for smooth optimization under stronger hypotheses, including linear convergence of gradient descent and local superlinear convergence of quasi-Newton methods \cite{2004_Boyd_BOOK,2006_Nocedal_BOOK}.
Geometric step rules can also yield linear convergence for subgradient methods under suitable conditions on the objective and the steps \cite{1977_Goffin,1985_Shor_BOOK}.

In this paper, we propose choosing the specular derivative $f^{\sd}(x_k)$ as the subgradient $g_k$.
Specular differentiation on $\mathbb{R}$ was introduced in \cite{2022a_Jung} and further developed in \cite{2026a_Jung}.
For a real-valued function $f:I\to\mathbb{R}$ on an open set $I\subset\mathbb{R}$, the specular derivative at $x\in I$ is given by
\begin{equation}\label{def:sd}
  f^{\sd}(x)
  :=
  \tan\left(
    \frac{\arctan\bigl(f'_+(x)\bigr)+\arctan\bigl(f'_-(x)\bigr)}{2}
  \right),
\end{equation}
whenever both one-sided derivatives $f'_+(x)$ and $f'_-(x)$ exist and are finite.
For the general definition, see \cite{2026a_Jung}.
In particular, $f^{\sd}(x)=f'(x)$ at every point of differentiability.
In the convex setting considered below, the representation \eqref{def:sd} always applies.

In particular, the specular derivative of a convex function is a subgradient at every point \cite[Cor.~2.21]{2026a_Jung}.
We therefore propose the \emph{specular gradient} (SPEG) method for \eqref{opt:1d}, which generates a sequence from a given initial point $x_0 \in \mathbb{R}$ by
\begin{equation}
  \label{eq:speg_update}
  x_{k+1} = x_k - \gamma_k f^{\sd}(x_k), \qquad k = 0, 1, 2, \ldots,
\end{equation}
where $f^{\sd}(x_k)$ denotes the specular derivative of $f$ at $x_k$ and $\gamma_k > 0$ is the step size.
This choice specifies the subgradient directly from the two one-sided derivatives, and hence no additional selection from $\partial f(x_k)$ is required at nondifferentiable points.
The resulting method retains the classical gradient update wherever $f$ is differentiable.

By combining the specular subgradient selection with normalized directions and geometrically decreasing step lengths, we establish R-linear convergence of the SPEG iterates, provided that the minimum is attained and the initial step length satisfies a suitable distance bound.
This yields an explicit subgradient method with a provable R-linear bound on the objective error.
In this one-dimensional setting, neither strong convexity nor differentiability of $f$ is required.
We use the terminology of \cite{2006_Nocedal_BOOK} for convergence rates.

The remainder of this paper is organized as follows.
\Cref{sec:convex_analysis_w_sd} presents convexity characterizations and a sufficient condition for minimality in terms of specular derivatives.
\Cref{sec:convergence} establishes R-linear convergence of SPEG with explicit error bounds for the iterates and objective values.
\Cref{sec:numerical_examples} describes the implementation and presents numerical comparisons on two convex examples.
Finally, \cref{sec:conclusion} discusses the scope of the results and possible extensions to higher dimensions.

\section{Convex analysis with specular derivatives}
\label{sec:convex_analysis_w_sd}

Throughout this section, let $I$ be an open interval in $\mathbb{R}$.

\begin{theorem}
  \label{thm:convexity}
  Let $f:I\to\mathbb{R}$ be continuous.
  The following are equivalent.
  \begin{enumerate}[label=\upshape(\alph*)]
    \item
    \label{thm:convexity-1}
    $f$ is convex on $I$.

    \item
    \label{thm:convexity-2}
    $f^{\sd}$ exists everywhere on $I$ and is nondecreasing.

    \item
    \label{thm:convexity-3}
    $f^{\sd}$ exists everywhere on $I$ and
    \begin{equation}\label{ineq: under graph}
      f(x)\geq f(y)+f^{\sd}(y)(x-y)
      \qquad\text{for all }x,y\in I.
    \end{equation}
  \end{enumerate}
\end{theorem}

\begin{proof}
  Assume \cref{thm:convexity-1}.
  By \cite[Cor.~2.21]{2026a_Jung}, $f$ is specularly differentiable on $I$ and $f^{\sd}(x)\in\partial f(x)$ for every $x\in I$.
  The one-dimensional characterization $\partial f(x)=[f'_-(x),f'_+(x)]$ therefore gives
  \begin{equation}\label{eq:sd_subgradient}
    f'_-(x)\leq f^{\sd}(x)\leq f'_+(x),
    \qquad f^{\sd}(x)\in\partial f(x)
    \quad(x\in I).
  \end{equation}
  For $x<y$, the inequality $f'_+(x)\leq f'_-(y)$ for convex functions yields
  \begin{displaymath}
    f^{\sd}(x)\leq f'_+(x)\leq f'_-(y)\leq f^{\sd}(y).
  \end{displaymath}
  This proves \cref{thm:convexity-2}, while the subgradient inequality gives \cref{thm:convexity-3}.

  Conversely, assume \cref{thm:convexity-2}.
  For $c<x<d$ in $I$, apply the Quasi-Mean Value Theorem \cite[Thm.~3.3]{2026a_Jung} on $[c,x]$ and $[x,d]$ to obtain points $u\in(c,x)$ and $v\in(x,d)$.
  Since $u<v$ and $f^{\sd}$ is nondecreasing, we have
  \begin{displaymath}
    \frac{f(x)-f(c)}{x-c}\leq f^{\sd}(u)
    \leq f^{\sd}(v)\leq\frac{f(d)-f(x)}{d-x}.
  \end{displaymath}
  Rearranging yields
  \begin{displaymath}
    f(x)\leq\frac{d-x}{d-c}f(c)+\frac{x-c}{d-c}f(d),
  \end{displaymath}
  which is convexity.

  Finally, assume \cref{thm:convexity-3}.
  Put $w=\lambda u+(1-\lambda)v$ for $u,v\in I$ and $0<\lambda<1$.
  Apply \eqref{ineq: under graph} with $y=w$ and $x=u,v$, respectively.
  Multiplication by $\lambda$ and $1-\lambda$ and addition give
  \begin{displaymath}
    \lambda f(u)+(1-\lambda)f(v)
    \geq f(w)+f^{\sd}(w)\bigl(\lambda(u-w)+(1-\lambda)(v-w)\bigr)
    =f(w).
  \end{displaymath}
  Thus $f$ is convex.
\end{proof}

\begin{corollary}[A sufficient condition for minimality]\label{crl:spd_is_zero}
  If $f:I\to\mathbb{R}$ is convex and $f^{\sd}(x^{\ast})=0$, then $x^{\ast}$ is a global minimizer of $f$ on $I$.
\end{corollary}

\begin{proof}
  Substitute $y=x^{\ast}$ into \eqref{ineq: under graph}.
\end{proof}

The converse need not hold.
For example,
\begin{equation}\label{ex : converse of spd is zero}
  f(x)
  =
  \begin{cases}
    x,  & x \geq 0,\\
    x^2,& x < 0
  \end{cases}
\end{equation}
is convex and has the unique minimizer $0$, but $f^{\sd}(0) = \sqrt{2} - 1$.
In practical implementations, one may terminate a gradient or subgradient method when the magnitude of the chosen gradient or subgradient falls below a prescribed tolerance.
However, the example above shows that the specular derivative need not vanish at a minimizer.
Thus, the stopping criterion $|f^{\sd}(x_k)| \leq \varepsilon$ may fail to recognize an exact minimizer.

\section{Convergence analysis}
\label{sec:convergence}

For a convex function $f:\mathbb{R}\to\mathbb{R}$ attaining its minimum, set
\begin{displaymath}
  f^{\ast} := \min_{x\in\mathbb{R}}f(x) 
  \qquad\text{and}\qquad
  \qquad X^{\ast} := \{x\in\mathbb{R}:f(x)=f^{\ast}\}.
\end{displaymath}
The minimizer set $X^{\ast}$ is nonempty, closed, and convex.
We write 
\begin{displaymath}
  \operatorname{dist}(x,X^{\ast})=\inf_{z\in X^{\ast}}|x-z|.
\end{displaymath}

SPEG need not be a descent method, as is also the case for general subgradient methods \cite{2014_Boyd}.
Indeed, for the function in \eqref{ex : converse of spd is zero}, a step from the minimizer $0$ leads to a negative point with a strictly larger objective value.
We therefore keep track of the best objective value obtained so far:
\begin{equation}\label{eq:updated_minimizer}
  f_k^{\mathrm{best}}:=\min_{0\leq j\leq k}f(x_j).
\end{equation}
The sequence $\{f_k^{\mathrm{best}}\}$ is nonincreasing and bounded below by $f^{\ast}$.

\begin{lemma}[Orientation outside the minimizer set]\label{lem : core inequality pointwise}
  Let $f:\mathbb{R}\to\mathbb{R}$ be convex with $X^{\ast}\neq\emptyset$.
  For every $x\in\mathbb{R}$, $z\in X^{\ast}$, and $g\in\partial f(x)$,
  \begin{equation}\label{eq:orientation}
    g(x-z)\geq f(x)-f^{\ast}\geq0.
  \end{equation}
  If $x\notin X^{\ast}$, then the first product is strictly positive.
  In particular, $g\neq0$ and $g$ has the sign of $x-z$.
\end{lemma}
\begin{proof}
  The subgradient inequality at $x$ gives $f(z)\geq f(x)+g(z-x)$.
  Rearrangement proves \eqref{eq:orientation}, and $f(x)>f^{\ast}$ when $x\notin X^{\ast}$.
\end{proof}

\begin{theorem}[R-linear convergence of SPEG]\label{thm:convergence_of_SPEG}
  Let $f:\mathbb{R}\to\mathbb{R}$ be convex with $X^{\ast}\neq\emptyset$.
  Choose $x_0\in\mathbb{R}$ and $t_0>0$ such that
  \begin{equation}\label{eq:initial_budget}
    \operatorname{dist}(x_0, X^{\ast}) \leq 2t_0 
    \qquad\text{and}\qquad
    \qquad t_k = t_0\left(\frac{1}{2}\right)^k \quad (k \geq 0).
  \end{equation}
  Generate the iterates by
  \begin{equation}\label{eq:normalized_speg}
    x_{k+1}=\begin{cases}
      x_k-t_k\displaystyle\frac{f^{\sd}(x_k)}{|f^{\sd}(x_k)|},& f^{\sd}(x_k)\neq0,\\[0.2cm]
      x_k,& f^{\sd}(x_k)=0.
    \end{cases}
  \end{equation}
  Then there exists $x^{\ast}\in X^{\ast}$ such that $x_k\to x^{\ast}$ and
  \begin{equation}\label{eq:speg_error_bound}
    \operatorname{dist}(x_k,X^{\ast})
    \leq
    |x_k-x^{\ast}|
    \leq
    2t_k
    =
    2t_0\left(\frac{1}{2}\right)^k
    \qquad(k\geq0).
  \end{equation}
  Moreover, there is a finite constant $L\geq0$ such that
  \begin{equation}\label{eq:objective_error}
    0
    \leq
    f_k^{\mathrm{best}}-f^{\ast}
    \leq
    f(x_k)-f^{\ast}
    \leq
    2Lt_0\left(\frac{1}{2}\right)^k
    \qquad(k\geq0).
  \end{equation}
  In particular, the iterates converge R-linearly to $x^{\ast}$.
  If $f^{\sd}(x_k)=0$ at some iteration, all subsequent iterates equal that minimizer.
\end{theorem}

\begin{proof}
  Put $d_k=\operatorname{dist}(x_k,X^{\ast})$.
  If $x_k\notin X^{\ast}$, let $z_k$ be its nearest point in the nonempty closed convex set $X^{\ast}$.
  By \eqref{eq:sd_subgradient} and \cref{lem : core inequality pointwise}, the step is directed toward $z_k$, and hence
  \begin{displaymath}
    d_{k+1}\leq |x_{k+1}-z_k|=|d_k-t_k|.
  \end{displaymath}
  If $x_k\in X^{\ast}$, then $d_{k+1}\leq|x_{k+1}-x_k|\leq t_k$.
  Therefore $d_k\leq2t_k$ implies $d_{k+1}\leq t_k=2t_{k+1}$ in both cases.
  Induction using \eqref{eq:initial_budget} gives $d_k\leq2t_k$ for every $k$.

  Moreover, for $m>k$,
  \begin{displaymath}
    |x_m-x_k|\leq\sum_{j=k}^{m-1}|x_{j+1}-x_j|
    \leq\sum_{j=k}^{m-1}t_j\leq2t_k.
  \end{displaymath}
  Thus $\{x_k\}$ is Cauchy and converges to some $x^{\ast}\in\mathbb{R}$.
  Since $d_k\to0$ and $X^{\ast}$ is closed, $x^{\ast}\in X^{\ast}$.
  Letting $m\to\infty$ gives \eqref{eq:speg_error_bound}.

  The explicit bound in \eqref{eq:speg_error_bound} establishes R-linear convergence of the iterates.
  All iterates and their limit lie in $K=[x_0-2t_0,x_0+2t_0]$ because the total step length is at most $2t_0$.
  Since $f$ is finite and convex on $\mathbb{R}$, it is Lipschitz on $K$ with some finite constant $L\geq0$; see \cite[Chap.~1]{1984_VanTiel_BOOK}.
  Hence, \eqref{eq:speg_error_bound} gives
  \begin{displaymath}
    0\leq f(x_k)-f^{\ast}
    \leq L|x_k-x^{\ast}|
    \leq2Lt_0\left(\frac{1}{2}\right)^k.
  \end{displaymath}
  Together with \eqref{eq:updated_minimizer}, this proves \eqref{eq:objective_error}.
  Finally, \cref{crl:spd_is_zero} proves the assertion about a zero specular derivative.
\end{proof}

\begin{remark}[Initialization from a known interval]\label{rmk:initialization}
  If a known interval $[a,b]$ with $a<b$ contains a minimizer, a natural choice is $x_0=\frac{1}{2}(a+b)$ and $t_0 = \frac{1}{4}(b - a)$.
  Then \eqref{eq:initial_budget} holds and all iterates remain in $[a,b]$.
\end{remark}

\begin{remark}[Necessity of the initial distance bound]\label{rmk:necessary_distance}
  The distance condition in \eqref{eq:initial_budget} is necessary for convergence to a minimizer under the prescribed step lengths.
  Indeed, the total step length is $\sum_{k\geq0}t_k=2t_0$, so every iterate and its limit lie within $2t_0$ of $x_0$; hence no minimizer can be reached if $\operatorname{dist}(x_0,X^{\ast})>2t_0$.
  For example, let $f(x)=|x-3|$, $x_0=0$, and $t_0=1$.
  Then the iterates satisfy $x_k=2(1-2^{-k})$ and converge to $2$, although the unique minimizer is $3$.
\end{remark}

\section{Numerical results}
\label{sec:numerical_examples}

For the numerical experiments, we use the SPEG implementation in \cite{py_1.3.2} with the geometric update rule in \cref{alg:SPEG}.
The code for reproducing these experiments is publicly available among the examples in the same repository.
The distance stopping criterion in this algorithm uses the bound established in \cref{thm:convergence_of_SPEG}.
The algorithm returns the current iterate, to which this distance bound applies, together with the best objective value observed.

For the fixed-budget comparison, we disable the distance stopping test in \cref{alg:SPEG} and keep the sequence constant after early termination.
The specular derivative is evaluated from analytic one-sided derivatives using a formula equivalent to \eqref{def:sd}.
All computations use double precision.

\begin{algorithm}
  \caption{\textsc{SPEG with geometric step lengths}}\label{alg:SPEG}
  \begin{algorithmic}[1]
    \REQUIRE{Convex $f:\mathbb{R}\to\mathbb{R}$ attaining its minimum, $x_0\in\mathbb{R}$, $t_0>0$ with $\operatorname{dist}(x_0,X^{\ast})\leq2t_0$, tolerance $\varepsilon>0$, maximum updates $N$}
    \ENSURE{Current iterate $x$, best value $f^{\mathrm{best}}$, distance bound $r$}
    \STATE{$x\gets x_0$, $t\gets t_0$, $k\gets0$, $f^{\mathrm{best}}\gets f(x_0)$}
    \WHILE{$k<N$ and $2t>\varepsilon$}
      \STATE{Compute $g\gets f^{\sd}(x)$}
      \IF{$g=0$}
        \RETURN{$x$, $f^{\mathrm{best}}$, $0$}
      \ENDIF
      \STATE{$\displaystyle x\gets x-t\frac{g}{|g|}$}
      \STATE{$\displaystyle t\gets\frac{t}{2}$, $k\gets k+1$}
      \STATE{$f^{\mathrm{best}}\gets\min\{f^{\mathrm{best}},f(x)\}$}
    \ENDWHILE
    \RETURN{$x$, $f^{\mathrm{best}}$, $2t$}
  \end{algorithmic}
\end{algorithm}
\FloatBarrier

For comparison, we use the \texttt{SGD} and \texttt{Adam} implementations in PyTorch \cite{2024_PyTorch} for gradient descent (GD) and Adam, respectively.
Adam \cite{2017_Kingma} is an adaptive stochastic gradient method.
GD uses no momentum, and Adam uses its default moment parameters $(0.9,0.999)$ and denominator constant $10^{-8}$.
Here both methods use automatic differentiation of the full objective, with the convention $\operatorname{sgn}(0)=0$ for the absolute-value terms.
Their learning rates have the form $\alpha_k=\alpha_0(k+1)^{-p}$ with $p\in\{0,\frac{1}{2},1\}$.
For each objective, we select $(\alpha_0,p)$ separately for GD and Adam from a finite grid by minimizing the mean objective error at $k=N$ over a separate set of $100$ initial points, sampled independently and uniformly from the same interval as the comparison points.

\begin{example}[Elastic Net]\label{ex:elastic_net}
  The first example is the Elastic Net, introduced in \cite{2005_Zou}.
  We consider its one-dimensional objective function $E:\mathbb{R}\to\mathbb{R}$, defined by
  \begin{equation}\label{ex:Elastic_Net}
    E(x)
    =\frac{1}{2m}\left\|Ax-\mathbf{b}\right\|_{\mathbb{R}^m}^{2}
    +\frac{\lambda_2}{2}x^2+\lambda_1|x|,
  \end{equation}
  where $m\geq1$, $A=(a_1,\ldots,a_m)^{\top}\in\mathbb{R}^{m\times1}$, $\mathbf{b}=(b_1,\ldots,b_m)^{\top}\in\mathbb{R}^{m}$, and $\lambda_1,\lambda_2\geq0$ are given.
  The optimization variable is scalar, while $m$ is the number of observations.
  If $\|A\|_{\mathbb{R}^m}^2+m\lambda_2>0$, its unique minimizer is
  \begin{displaymath}
    x_E^{\ast}
    =\frac{\operatorname{sgn}(A^{\top}\mathbf{b})\max\{|A^{\top}\mathbf{b}|-m\lambda_1,0\}}{\|A\|_{\mathbb{R}^m}^2+m\lambda_2}.
  \end{displaymath}
  For the reported experiment, we set $m=100$, draw $a_1,\ldots,a_m$ independently from the standard normal distribution, and take $b_i=1$ for every $i$.
  This single realization of $A$ and $\mathbf{b}$ is shared by all methods and initial points.
  We choose $\lambda_1=1$ and $\lambda_2=\frac{1}{2}$, which give $x_E^{\ast}=0$ and $E^{\ast}=\frac{1}{2}$ for this realization.

  All three methods use the same $100$ initial points sampled uniformly from $[-4,4]$, with $N=1000$ updates.
  SPEG uses $t_k=t_0 2^{-k}$ with $t_0=3$, and hence \eqref{eq:initial_budget} holds for every initial point in $[-4, 4]$.
  The selected learning rates are $\frac{3}{10(k+1)}$ for GD and $\frac{2}{3(k+1)}$ for Adam.

  \begin{figure}[htbp]
    \centering
    \includegraphics[width=\textwidth]{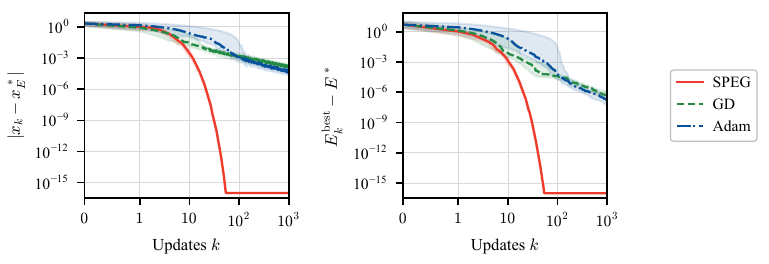}
    \par\vspace{-10pt}
    \caption{Elastic Net \eqref{ex:Elastic_Net} with $m=100$, $\lambda_1=1$, and $\lambda_2=\frac{1}{2}$. Left: current-iterate error. Right: best objective error. Curves and bands show medians and interquartile ranges over $100$ initial points for the fixed random design.}
    \label{fig:elastic_net_results}
  \end{figure}

  \begin{table}[htbp]
    \centering
    \caption{Median errors for the Elastic Net at $k=1000$ over $100$ initial points.}
    \label{tab:elastic_net_results}
    \small
    {\begin{tabular}{lcc}
\toprule
Method & Median $|x_k-x_{E}^\ast|$ & Median $E_k^{\mathrm{best}}-E^\ast$ \\
\midrule
SPEG & $1.87\times10^{-301}$ & $1.50\times10^{-301}$ \\
GD & $1.29\times10^{-4}$ & $5.48\times10^{-7}$ \\
Adam & $4.44\times10^{-5}$ & $1.76\times10^{-7}$ \\
\bottomrule
\end{tabular}
}
  \end{table}
\end{example}

\begin{example}[Sum of absolute values]\label{ex:absolute_sum}
  The second objective function $F:\mathbb{R}\to\mathbb{R}$ is defined by
  \begin{equation}\label{eq:absolute_sum}
    F(x)=\sum_{i=0}^{99}\left(
      \left|x-\frac{i}{100}\right|+\left|x+\frac{i}{100}\right|
    \right).
  \end{equation}
  This function is convex and piecewise linear, with $199$ nondifferentiable points, and is not strongly convex.
  Since $F(x) \geq 99 + 2|x|$ and $F(0) = 99$, its unique minimizer is $x_F^{\ast} = 0$ and its minimum is $F^{\ast} = 99$.

  All three methods use the same $100$ initial points sampled uniformly from $[-1,1]$, with $N=50$ updates.
  SPEG uses $t_k=t_0 2^{-k}$ with $t_0=\frac{1}{2}$, and hence \eqref{eq:initial_budget} holds for every initial point in $[-1,1]$.
  The selected learning rates are $\frac{1}{200(k+1)}$ for GD and $\frac{3}{10(k+1)}$ for Adam.

  \begin{figure}[htbp]
    \centering
    \includegraphics[width=\textwidth]{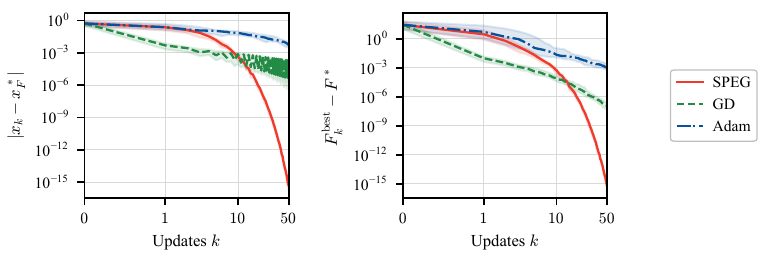}
    \par\vspace{-10pt}
    \caption{Sum of absolute values \eqref{eq:absolute_sum}. Left: current-iterate error. Right: best objective error. Curves and bands show medians and interquartile ranges over $100$ initial points.}
    \label{fig:absolute_sum_results}
  \end{figure}

  \begin{table}[htbp]
    \centering
    \caption{Median errors for the sum of absolute values at $k=50$ over $100$ initial points.}
    \label{tab:absolute_sum_results}
    \small
    {\begin{tabular}{lcc}
\toprule
Method & Median $|x_k-x_{F}^\ast|$ & Median $F_k^{\mathrm{best}}-F^\ast$ \\
\midrule
SPEG & $4.44\times10^{-16}$ & $4.44\times10^{-16}$ \\
GD & $1.96\times10^{-4}$ & $1.00\times10^{-7}$ \\
Adam & $5.49\times10^{-3}$ & $1.01\times10^{-3}$ \\
\bottomrule
\end{tabular}
}
  \end{table}
\end{example}

\Cref{fig:elastic_net_results,fig:absolute_sum_results} show the errors $|x_k-x_E^{\ast}|$ and $E_k^{\mathrm{best}}-E^{\ast}$ for $E$, and the corresponding errors for $F$.
Here $E_k^{\mathrm{best}}$ and $F_k^{\mathrm{best}}$ are defined as in \eqref{eq:updated_minimizer}.
Curves show the median over $100$ runs, and shaded bands span the 25th to 75th percentiles.
Objective errors are evaluated using algebraically equivalent formulas that avoid subtracting nearly equal function values.
For plotting, values below $10^{-16}$ are displayed at $10^{-16}$.
The SPEG entries in \cref{tab:elastic_net_results,tab:absolute_sum_results} are compatible with double precision because the minimizers are exactly $x_E^{\ast}=x_F^{\ast}=0$ and the objective errors are computed from the stable formulas above; the current-iterate errors agree with the bounds $2t_{1000} = 6 \cdot 2^{-1000} \approx 5.6\times 10^{-301}$ and $2t_{50} = 2^{-50} \approx 8.9\times 10^{-16}$ from \cref{thm:convergence_of_SPEG}.

At $k=1000$ for $E$ and $k=50$ for $F$, SPEG has smaller median current-iterate and best objective errors than GD and Adam; see \cref{tab:elastic_net_results,tab:absolute_sum_results}.
The initial decrease of the SPEG errors is consistent with the R-linear estimates in \cref{thm:convergence_of_SPEG}.

\section{Scope and further directions}
\label{sec:conclusion}

We establish R-linear convergence of SPEG for one-dimensional convex optimization using normalized steps of geometrically decreasing length.
The result assumes that the minimum is attained and the initial distance bound \eqref{eq:initial_budget} holds.
It allows nonunique minimizers and requires neither differentiability nor strong convexity.
Numerical experiments on the Elastic Net and a sum of absolute values illustrate the convergence behavior.
Extending this analysis to higher dimensions remains a direction for future work, since the one-dimensional orientation argument does not directly apply there.

\bibliographystyle{siamplain}
\bibliography{reference}

\end{document}